\documentclass{article}

\usepackage{graphicx,fullpage}%
\usepackage{multirow}%
\usepackage{amsmath,amssymb,amsfonts}%
\usepackage{amsthm}%
\usepackage{mathrsfs}%
\usepackage[title]{appendix}%
\usepackage{xcolor}%
\usepackage{textcomp}%
\usepackage{manyfoot}%
\usepackage{booktabs}%
\usepackage{algorithm}%
\usepackage{algorithmicx}%
\usepackage{algpseudocode}%
\usepackage{listings}%
\usepackage[T1]{fontenc}
\usepackage[utf8]{inputenc}
\usepackage{lmodern}
\usepackage{soul}  
\usepackage{multirow}   
\usepackage{makecell}

\usepackage{cancel}
\usepackage{authblk}
\usepackage{url}

\providecommand{\DIFdel}[1]{}
\providecommand{\DIFdel}[1]{\ignorespaces}

\DeclareMathOperator*{\argmin}{\arg\min}

\theoremstyle{thmstyleone}%
\theoremstyle{thmstyletwo}%

\theoremstyle{thmstylethree}%

\begin{document}

\title{A Differentiable Simulation of the Eye for Patient-Specific Strabismus Surgery Planning}


\author[1]{Even Harsigny}
\author[1]{Pablo Alvarez}
\author[1]{Michel Duprez}
\author[1]{Stéphane Cotin}
\affil[1]{Université de Strasbourg, CNRS, Inria, ICube, F-67000 Strasbourg, France}




\maketitle

\abstract{
\textbf{Purpose:} Up to 4\% of adults will develop strabismus in their lifetime. The most common surgical intervention involves adjusting the length of one or more extraocular muscles to correct the angular deviation. This correction depends on surgical expertise and statistical reference tables, which often fail to yield optimal results for patients with atypical eye morphology.
Our work proposes a physics-based modeling approach to personalized surgical planning, accounting for patient-specific eye anatomy. 

\textbf{Methods:} We built a physics-based simulator of the eye and its muscles, incorporating patient-specific geometry and Hill-type muscle biomechanics. We solve an optimization problem to find the surgical dosage that minimizes angular deviation. The model is implemented as a fully differentiable simulation, enabling efficient optimization. We validated the framework by comparing its predictions with standard surgical tables for emmetropic eyes before applying it to anatomically atypical virtual patients.

\textbf{Results:} Our model's predictions for emmetropic eyes were first validated, demonstrating a strong fit with standard surgical tables. More importantly, for high-myopia models, the framework computed  a clinically significant increase in the required surgical dosage compared to standard eyes. This computed recession difference is highly relevant as surgical plans are adjusted in 0.5~mm increments. 

\textbf{Conclusion:} Our results show that our model provides a calibrated surgical plan that, unlike standard tables, also accounts for pathologies involving atypical eye shapes. This patient-specific model represents a step toward personalized surgical planning, with the potential to improve dosage accuracy and surgical outcomes for atypical cases.
}

\section{Introduction}\label{sec:introduction}

Strabismus is a disorder in which the eyes do not align when looking directly at an object. This pathology requires surgical intervention when optical correction is insufficient to restore binocular alignment and eliminate diplopia. The primary goal of surgery is to adjust extraocular muscle (EOM) tension, typically via recession, to correct the deviation. However, high rates of re-operation due to under- or over-correction are common, as outcomes remain unpredictable \cite{Niu2024}. While perioperative techniques such as adjustable sutures \cite{Zhang2012} enable postoperative refinement to mitigate surgical variability, they cannot fully compensate for fundamental planning errors caused by atypical anatomy. Accuracy is critical, as residual deviations as small as two prism diopters (PD) can cause persistent diplopia \cite{Wright2009}. Yet variability in surgical outcomes often exceeds this threshold because of a lack of consensus on optimal target criteria \cite{Hatt2012}.

When it comes to planning strabismus surgery, surgeons strongly rely on their experience and on nomograms. However, as strabismus surgery is individualized, guidelines cannot be applied to every case \cite{Khurana2018}; instead, they serve as references. Indeed, the statistical values reported in tables, as in \cite{Khurana2018}, apply to average eye anatomies and do not account for patient-specific variations in muscle strength, insertion points, or overall globe geometry. This approach tends to be problematic for patients with atypical anatomy.
A simulator of the eye would enable surgeons to plan their procedures virtually, while accounting for patient-specific characteristics. 

Oculomotor modeling literature presents several approaches for surgical planning, which can be grouped into three main biomechanical model families by complexity.

First, static geometric models established the foundation of computerized surgical planning. Initiated by Robinson \cite{Robinson1975}, who solved for static force equilibrium in the orbit, these approaches evolved into more sophisticated formulations, such as the SQUINT model \cite{Miller1984}. These works led to the development of frequently used clinical software tools such as \textit{Orbit} and \textit{SEE++} \cite{Haslwanter2005}, which allow surgeons to simulate surgical outcomes virtually. However, while these platforms provide a graphical interface for planning, they are inherently limited to static fixations and offer only forward prediction based on user inputs.

Second, Finite Element (FEM) and strand-based models \cite{Wei2010} provide high bio-fidelity by simulating tissue deformation, but are too computationally intensive for iterative, patient-specific planning.

Third, multi-body dynamic models \cite{Iskander2018a} represent a balanced alternative, modeling the eye as a rigid body actuated by Hill-type muscles and incorporating active pulley constraints. However, they have primarily been used for forward kinematic analysis rather than automated surgical planning, as they lack the differentiable formulation needed to efficiently solve the inverse problem of finding optimal surgical parameters.

Additionally, data-driven methods (statistical \cite{Leite2021} and AI \cite{Mao2021}) predict surgical plans from patient data patterns. While effective for typical cases, these approaches rely on extensive training data and may be unreliable for rare anatomies such as high myopia with strabismus, where data are limited.

Specific studies have also investigated the impact of anatomical variations on surgical planning to determine if standard tables could be adjusted using simple biometric factors. Hirnschall et al. \cite{Hirnschall2022} and Beisse et al. \cite{Beisse2020} analyzed the correlation between axial length and surgical dosage. They derived linear adjustment rules (or noted that optical and anatomical factors partially cancel each other out \cite{Beisse2020}), suggesting that, for the general population, the eye acts as a uniformly scaled sphere, with biometric scaling playing a relatively minor role.

However, the linear assumption breaks down in pathological cases like high myopia. Hoerantner et al. \cite{Hoerantner2007} showed that high myopic eyes may exhibit geometric distortions, not simple scaling, which alter muscle paths. Thus, while linear scaling may suit mild variations, substantial anatomical changes in high myopes likely necessitate physics-based force equilibrium modeling, as investigated here.

We present a novel, differentiable physics-based framework for patient-specific strabismus surgery planning that, unlike static geometric or computationally intensive finite element models, enables efficient, gradient-based optimization while maintaining biomechanical fidelity. Our main contributions are threefold:
\begin{enumerate}
  \item A differentiable simulator of the eye and its muscles integrating active pulley mechanics with 3D globe dynamics.
  \item A three-step inverse problem methodology: first, to compute the muscular activations in primary gaze for a healthy eye, second, to identify the patient-specific pathology (muscle weakness) from a clinical measurement (prism diopters), and third, to predict the optimal recession distance.
  \item A demonstration of clinical relevance by quantifying the impact of high myopia on surgical planning. Our experiments reveal that standard surgical tables lead to systematic under-correction in these cases.
\end{enumerate}

\section{Methods}\label{sec:methods}

\subsection{Biomechanical and anatomical model}

\subsubsection{Ocular globe}
We model the eye as a rigid-body ellipsoid defined by its three radii (axial, height, and width). These radii can be personalized to adapt the model to specific patient anatomies. When all three radii are equal, the model simplifies to a sphere, which we use as our baseline healthy eye. In particular, this baseline eye has a radius of 12~mm, which is approximately the mean radius of healthy eyes \cite{Augusteyn2012}. We constrain the globe at its center of mass to enable rotations only, with three rotational degrees of freedom. We model the eye as a rigid body, consistent with established dynamic frameworks \cite{Wei2010}.
\begin{figure}[h!]
    \centering
    \includegraphics[width=1.0\textwidth]{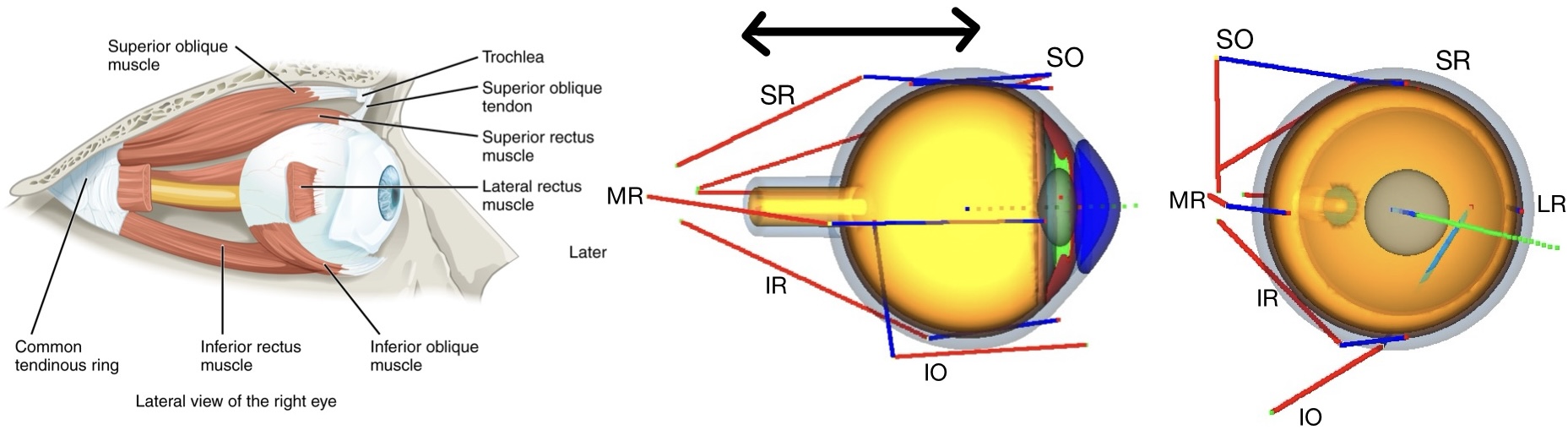}
    \caption{(Left) Anatomical reference of the extraocular muscles (EOM), source: \texttt{https://commons.wikimedia.org/wiki/File:1412\_Extraocular\_Muscles.jpg}. (Center, Right) Side and frontal views of our biomechanical model implemented in SOFA, showing the 12~mm radius emmetropic globe. The simplified muscle paths are clearly shown, with the red lines representing the origin-to-pulley segments and the blue lines representing the pulley-to-insertion segments. The black arrow indicates the movement of pulleys}
    \label{fig:eom}
\end{figure}

\subsubsection{Extraocular muscles}

Eye movements are governed by six extraocular muscles, all incorporated in our model: the Medial Rectus (MR), Lateral Rectus (LR), Superior Rectus (SR), Inferior Rectus (IR), Superior Oblique (SO), and Inferior Oblique (IO) (Fig. \ref{fig:eom}). Each muscle originates at the common tendinous ring and inserts near the sclera, traversing connective tissue pulleys. Consistent with established modeling approaches \cite{Iskander2018a}, we represent the muscle path as two straight segments, as shown in Fig. \ref{fig:eom}. Muscle contraction produces torque, resulting in eye rotation. Pulleys act as dynamic origins, guiding muscle trajectories and ensuring correct force transmission. We implement the active pulley system, which permits coordinated pulley motion during eye rotation, following \cite{Miller1984} and \cite{Iskander2018a}.

To compute muscle force, we employ the Hill-type model \cite{Millard2013}, which includes passive $f_{P}$ (spring-like, length-dependent) and active $f_{A}$ (activation and length-dependent, velocity-modulated) force components with activation level $a(t) \in [0,1]$. The tendon is modeled as rigid, consistent with prior studies \cite{Iskander2018a}. In our model, a muscle path is defined by three points: the anatomical (fixed) origin $\mathbf{p_{\text{orig}}}$, the active pulley position $\mathbf{p_{\text{pulley}}(x)}$, which is a function of gaze, and the position of the insertion $\mathbf{p_{\text{ins}}(x)}$, also a function of gaze, with $\mathbf{x} \in SO(3)$, as it rotates the globe. Thus, the length of each muscle is defined as: 
\begin{equation}
    l_M(\mathbf{x}; l_{\text{opt}}) = \frac{\| \mathbf{p}_{\text{pulley}}(\mathbf{x}) - \mathbf{p}_{\text{orig}} \| + \| \mathbf{p}_{\text{ins}}(\mathbf{x}) - \mathbf{p}_{\text{pulley}}(\mathbf{x}) \|}{l_{\text{opt}}},
    \label{eq:length}
\end{equation}
with $l_{\text{opt}}$ being the optimal length of a muscle at which it produces the maximum active force \cite{Miller1984}. To compute the full force generated by a muscle, the active and passive forces are multiplied by the maximum isometric force parameter $F_{\text{max}}$, which defines the maximal force a muscle can produce:
\begin{equation}
    F_{\text{total}} = F_{\text{max}} \left( f_{A}(l_M)f_{V}(v_M)a(t) + f_{P}(l_M)  \right).
\label{eq:hill}
\end{equation}

The active $f_{A}(l_M)$ and passive $f_{p}(l_M)$ force-length relationships were fitted to match the specific extraocular muscle properties reported by Miller and Robinson \cite{Miller1984}, ensuring the model captures the non-linearities specific to EOM. In this static study, the velocity $v_{M}$ is null, which gives $f_{V}(v_M)=1$. We digitized the original data points and implemented a 3rd-order spline interpolation to define the continuous passive and active force functions. This approach ensures an exact match with the physiological reference curves \cite{Miller1984}.

\subsubsection{Eyeball motion}

Here, we show how the global eye motion is computed from the elongation and contraction of the 6 EOMs. For a given muscle $i$ exerting a force $F_{\text{total}, i}$, the resulting torque $\tau_i$ acting on the eyeball is a function of gaze and activation: 
\begin{equation}
    \tau_i(\mathbf{x}, a_i; F_{\text{max}, i}) = \mathbf{r}_i(\mathbf{x}) \times \left( F_{\text{total}, i}(\mathbf{x}, a_i; F_{\text{max}, i}) \cdot \mathbf{d}_i(\mathbf{x}) \right),
\end{equation}
with $\mathbf{r}_i(\mathbf{x})$ the moment arm of the muscle (\textit{i.e.} the vector going from the center of the eye to the insertion point $\mathbf{p}_{\text{ins}, i}(\mathbf{x})$),  and $\mathbf{d}_i(\mathbf{x})$ the muscle direction vector from the insertion point to the pulley point. The net torque applied to the eye is the sum of all six torques, plus a passive torque due to the surrounding tissues resisting eye movement \cite{Robinson1975}: 
\begin{equation}
    \tau_{\text{net}}(\mathbf{x}, \mathbf{a}; \mathbf{F_{\text{max}}}) = \sum_{i=1}^{6} \tau_i(\mathbf{x}, a_i; F_{\text{max}, i}) +  \tau_{\text{passive}},
\label{eq:sum_torque}
\end{equation}
where $\mathbf{a}$ and $\mathbf{F_{\text{max}}}$ are the aggregation of muscle activations $a_i$ and exerting forces $F_{\text{max}, i}$, respectively. If muscle activations generate zero net torque for a given gaze, the eye remains in stable equilibrium. Non-zero net torque causes the eye to rotate accordingly. This framework explains how the baseline in a healthy, emmetropic eye adapts to anatomical variations, such as high myopia.

\subsubsection{Patient-specific adaptation}

Our model is designed to accommodate anatomical variations, with myopia serving as a key example due to its prevalence and significant increase in axial length. In myopia, the eye becomes more ellipsoid, with increases in axial length, height, and width quantified by Atchison et al. \cite{Atchison2004} as 0.35~mm/D, 0.19~mm/D, and 0.10~mm/D, respectively. This elongation alters both globe biomechanics and extraocular muscle (EOM) paths, increasing muscle forces and modifying moment arms, which impacts surgical planning. We modeled high myopia by non-uniformly scaling the globe, adjusting muscle insertions accordingly, while retaining pulley coordinates as in the emmetropic eye, since pulleys are orbitally suspended. This assumes axial elongation affects globe dimensions but not orbital connective tissues, barring specific pathologies.

\subsection{Optimization via differentiable simulation}

In this section, we address the numerical aspects of the optimization problem for predicting the muscle recession required to correct strabismus. We follow a three-step approach. The first step focuses on estimating the muscular activations that provide a zero net torque and fixate the eye in primary gaze (\textit{i.e.} zero eyeball rotation). The second step aims to model strabismus by determining the muscle parameters that use primary-gaze activations to produce a desired deviation. Finally, the third step is to find the optimal recession that corrects for this deviation. In this work, only exotropia and esotropia are studied, limiting deviations to the horizontal plane.

To solve these inverse problems, we implemented our simulation in the SOFA framework \cite{Faure2012}. While simpler solvers exist, SOFA was selected for its validated library of medical biomechanics constraints and its modular architecture. Furthermore, to enable gradient-based optimization, the entire analytic chain from the muscle's geometric points to the final net torque in Eq.~\eqref{eq:sum_torque} is written in \texttt{Python} using the JAX library within SOFA. This fully differentiable formulation allows the net torque equation to be differentiated with respect to any parameter, enabling efficient gradient computation to automatically find optimal surgical parameters, thus avoiding the trial-and-error approach of non-differentiable models. The following sections describe the specifics of each optimization problem.

\subsubsection{Muscular activation estimation}

The first step focuses on physiological stability involving the entire oculomotor system. We solve for the simultaneous muscular activations $a \in [0,1]^6$ of all six EOMs that result in zero net torque (Eq.~\eqref{eq:sum_torque}) to fixate the eye in a given gaze. To retrieve these optimal activations $\mathbf{a}^*$, the following inverse problem is solved:
\begin{equation}
\label{eq:fixation_problem}
   \mathbf{a}^* = \argmin_{\mathbf{a} \in [0, 1]^6}  \;\;\;  \left\| \tau_{\text{net}}(\mathbf{x}, \mathbf{a}; \mathbf{F_{\text{max}}}) \right\|^2 + \alpha \left\|\mathbf{a}\right\|^2,
\end{equation}
which minimizes the squared magnitude of the net torque at the desired gaze $\mathbf{x}$ while regularizing muscle activation energy. 

By solving \eqref{eq:fixation_problem} for $\mathbf{x} = \mathbf{0}$ on a healthy eyeball, we obtain the muscle activations required for holding primary gaze in a healthy subject, which we denote $\mathbf{a_{\text{healthy}}}$. In the following, we assume that, in the case of strabismus, the same muscle activations $\mathbf{a_{\text{healthy}}}$ no longer hold a primary gaze but rather a deviated gaze.

\subsubsection{Strabismus parameter optimization}

The goal of this second step is to simulate a virtual eye with a specific ocular deviation. Although our differentiable framework is generic and could optimize multiple parameters simultaneously (e.g., unilateral recession-resection), we deliberately constrain the optimization in this study to the maximum isometric force $F_{\text{max}}$ of a single horizontal muscle (either LR or MR). Clinically, this choice implies modeling the strabismus through its mechanical manifestation at the effector level, regardless of its neurological origin. Since the surgical procedure evaluated here (unilateral recession) relies on mechanically weakening a specific muscle to treat the deviation, representing the pathology as an increase in that muscle's force provides a direct biomechanical equivalent for surgical planning. 

For an eye with strabismus (one weaker or stronger $F_{\text{max}}$ for a horizontal EOM) with $\mathbf{a_{\text{healthy}}}$ muscle activations at the primary gaze, the developed net torque from Eq. \eqref{eq:sum_torque} is not zero, entailing a horizontal eye deviation. Let $D : \mathbf{(F_{\text{max}}}; \mathbf{a}) \mapsto \text{D}_{\text{stable}}$ be a function that maps a given set of maximal isometric forces $\mathbf{F_{\text{max}}}$ to a stable eye deviation $\text{D}_{\text{stable}}$ when a set of activations $\mathbf{a}$ are applied at primary gaze.  
Thus, for a healthy eye with healthy maximal isometric forces, the obtained stable eye deviation $\text{D}_{\text{stable}}$ would be zero, \textit{i.e.} the primary gaze. We model strabismus with a desired deviation D$_{\text{target}}$ by solving the following optimization problem: 
\begin{equation}
\label{eq:strabismus_problem}
    \mathbf{F^*_{\text{max}}} = \argmin_{\mathbf{F_{\text{max}}}} \left( D(\mathbf{F_{\text{max}}}; \mathbf{a} = \mathbf{a_{\text{healthy}}}) - \text{D}_{\text{target}} \right)^2.
\end{equation}

This process allows for graduating the degree of exotropia or esotropia by solving Eq.~\eqref{eq:strabismus_problem} for varying degrees of deviations D$_{\text{target}}$.

\subsubsection{Optimization of the recession amplitude}

Consistent with the single-muscle pathology modeled in the previous step, we determine the optimal recession for this specific muscle to return the eye to alignment. To determine this recession amplitude, a third inverse problem needs to be solved, and this is where the differentiable nature of our model is an advantage. Indeed, since the entire model is implemented on top of the JAX library, we can compute the gradient of our objective function ``for free'', reducing computational errors and improving optimization convergence compared to gradient-free methods.

This formulation enables the full optimization pipeline to converge in 18.3 seconds on a standard workstation. Although planning is not time-critical, such efficiency is essential for rapid clinical feedback and for enabling large-scale analyses, including the sensitivity analysis described in Section 3, which would be much more time-consuming with non-differentiable solvers.

The goal of the surgery is to stabilize the eye in primary gaze, which constitutes the standard definition of alignment success \cite{Hatt2012}. Consequently, we formulate the optimization problem in this primary position. The preservation of motility in secondary gazes (i.e., avoiding incomitance) is incorporated not as a simultaneous optimization target, but by verifying that the resulting recession remains within clinical safety limits (typically 10~mm \cite{Kassem2024}), as detailed in Section 3. Being in primary gaze means that the net torque generated by the EOMs must be 0 when the eye is in this position, while the muscles are activated with $\mathbf{a_{\text{healthy}}}$. Therefore, the objective is to find the new, optimal insertion point of the affected muscle that results in a net torque of zero, given $F^*_{\text{max}}$ as computed above.
From Eq.~\eqref{eq:length}, the new insertion point $\mathbf{p}_{\text{ins}}$ of the operated muscle becomes a parameter of the function $l_M$. Therefore, the function $ \tau_{\text{net}}$ now also depends on $\mathbf{p}_{\text{ins}}$, and we solve for the optimal insertion $\mathbf{p}^*_{\text{ins}}$ by minimizing the squared magnitude of the net torque at primary gaze :
\begin{equation}
    \mathbf{p}^*_{\text{ins}} = \argmin_{\mathbf{p}_{\text{ins}}} \left\| \tau_{\text{net}}(\mathbf{x}=\mathbf{0}, \mathbf{a} = \mathbf{a_{\text{healthy}}, \mathbf{p}_{\text{ins}}}; \mathbf{F_{\text{max}}} = \mathbf{F^*_{\text{max}}}) \right\|^2.
\end{equation}

Once $\mathbf{p}^*_{\text{ins}}$ is optimized, the geodesic distance between old and new insertion points is calculated, indicating the required surgical recession.

\section{Results and discussion}\label{sec:results}

We assessed agreement of the emmetropic (control) model with published clinical dosing tables and ranges \cite{Kassem2024, Wang2010, Lekskul2021}. First, we gathered minimal and maximal recession values from the literature to establish clinical recession intervals for each PD deviation. 
We then report consistency in terms of prediction-recession error relative to those intervals. Specifically, for each considered PD deviation, the error was quantified as the distance of the predicted recession to the nearest interval boundary. All predictions within the tabulated intervals were considered clinically acceptable and, therefore, yielded an error of zero.   
Within the common unilateral surgery range (15–25 PD) \cite{Kassem2024}, we achieved an average error of 0.15~mm, well below the typical 0.5~mm surgical increment.

For deviations above 25 PD, our model predicts a required recession exceeding standard surgical safety limits, which plateau at 10~mm \cite{Kassem2024}. 
For deviations exceeding 30 PD and under the unilateral recession constraint R $\leq$ 10~mm, our model predicts that full correction is unattainable. While this limit is imposed to prevent lateral incomitance (i.e., difference in horizontal deviation between lateral gaze), our model reveals that adhering to it inevitably compromises primary gaze alignment.
This observation is consistent with the under-correction rates reported for large-angle exotropia surgery \cite{Kassem2024}.
\begin{table}[htbp]
\centering
\caption{Comparison of projected surgical outcomes for a large-angle exotropia (35 PD) on the LR muscle. We compare the residual deviation obtained when applying the standard surgical table value (10~mm safety limit) versus the outcome obtained using the predicted recession}
\label{tab:failure_analysis}
\begin{tabular}{l l c c c}
\toprule
\textbf{Patient Model} & \textbf{Surgical Strategy} & \textbf{\thead{Applied \\ Recession}} & \textbf{\thead{Residual \\ Deviation}} & \textbf{Outcome} \\
\midrule
Emmetropic (Control) & Standard (safety limit) & 10.00~mm & 6.79 PD & Failure \\
Myopic ($-9\text{D}$) & Standard (safety limit) & 10.00~mm & 8.23 PD & Failure \\
\textbf{Myopic ($-9\text{D}$)} & \textbf{Patient-Specific (ours)} & \textbf{12.77~mm} & \textbf{0.00 PD} & \textbf{Success} \\
\bottomrule
\end{tabular}
\end{table}

To quantify the impact of these generic limitations on atypical anatomy, we simulated a large-angle exotropia (35 PD) scenario on the LR muscle (Table \ref{tab:failure_analysis}). Crucially, adhering to the 10~mm safety cap results in guaranteed surgical failure:
while an emmetropic eye is under-corrected (6.79 PD residual), the high myope ($-9\text{D}$) suffers a severe 8.23 PD residual deviation. Our model computes that a recession of 12.77~mm is geometrically necessary to correct the $-9\text{D}$ eye. 
Rather than advocating to exceed safe limits, this result serves as a model-based indicator that unilateral recession within the 10~mm cap may be insufficient, motivating alternative strategies (e.g., bilateral surgery). Consequently, our framework offers dual validation: it confirms standard surgical plans for common deviations in unilateral surgery ($\leq$ 25 PD) and reliably detects biomechanical impasses for atypical, large-angle cases.
\begin{figure}[h!]
    \centering
    \includegraphics[width=0.9\textwidth]{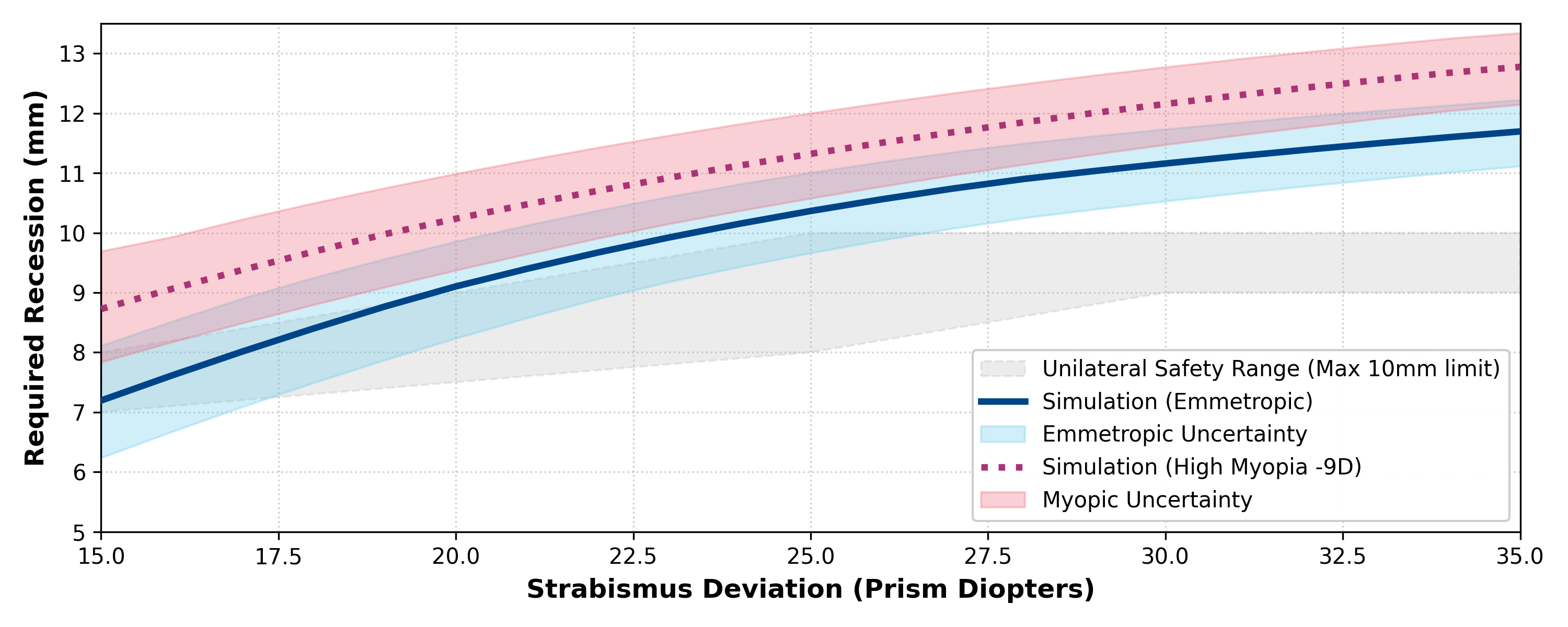}
    \caption{Comparison of simulated LR recession requirements for exotropia (blue: emmetropic, red: $-9\text{D}$ myopia) against standard surgical guidelines (gray zone). Shaded areas represent model uncertainty}
    \label{fig:chart_myopia}
\end{figure}

\begin{figure}[h!]
    \centering
    \includegraphics[width=0.9\textwidth]{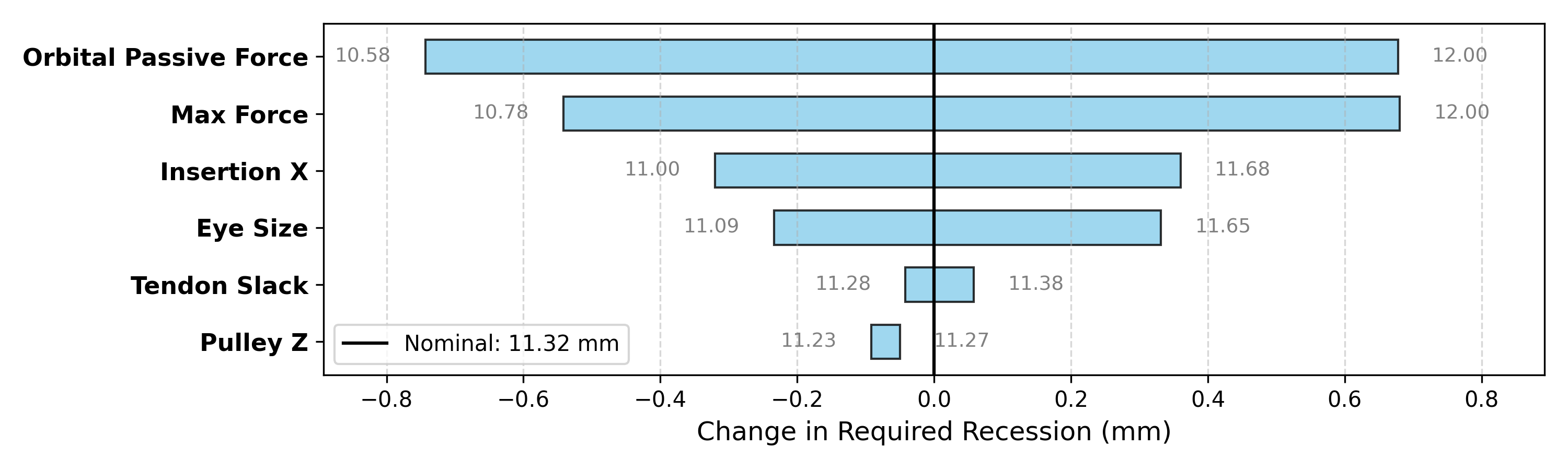}
    \caption{Sensitivity analysis (Tornado plot) for a $-9\text{D}$ myopic eye with a 25 PD deviation. The horizontal bars quantify the change in the required surgical recession (in mm) when each parameter is varied by ±1 standard deviation (SD) from its nominal value. The parameters are ranked by influence, identifying orbital passive force and maximum isometric force as the primary determinants of the surgical plan}
    \label{fig:tornado}
\end{figure}

To evaluate our model across axial lengths, we generated myopic patient models ($-3\text{D}$, $-6\text{D}$, $-9\text{D}$). Fig.~\ref{fig:chart_myopia} compares personalized plans for emmetropic and $-9\text{D}$ eyes, with the latter requiring 0.9--1.5~mm more recession, a clinically meaningful difference given standard 0.5~mm increments. To assess the robustness of this finding, we performed a univariate sensitivity analysis (N=1176) by perturbing the six parameters in Table~\ref{tab:sensitivity_params} individually to their upper and lower bounds ($\pm 1$ SD) while keeping others fixed to their nominal value. These variations generate the uncertainty bands shown in Fig.~\ref{fig:chart_myopia}, where, despite partial overlaps, the myopic model consistently requires larger recessions than the emmetropic baseline over the evaluated range. Finally, the relative influence of these parameters is ranked on the tornado diagram in Fig.~\ref{fig:tornado}. Results identify orbital passive force and maximum isometric force as the main determinants, regardless of eye length. The results in Fig.~\ref{fig:chart_myopia} only display recessions for emmetropic and $-9\text{D}$ myopic eyes for clarity. In contrast, the $-3\text{D}$ myopic eye showed 0.3–0.4~mm greater recession than the emmetropic eye, and the $-6\text{D}$ myopic eye required 0.6–1.0~mm more recession across the studied deviation angle range. This reveals a dual dependency: while eye geometry shifts the nominal surgical target (requiring more recession for myopes), the system's stability is governed by tissue mechanics. 
This supports the conclusion that simple geometric scaling alone may be insufficient, motivating physics-based modeling.

Furthermore, our simulation quantifies a biomechanical phenomenon that explains clinical variability: surgical efficacy drops by 17.3\% in small deviations (from 2.08 to 1.72 PD/mm) for high myopes. This aligns with Gezer et al. \cite{Gezer2004}, who reported that high myopes respond less favorably to surgery than emmetropes for the same deviation. This confirms that the discrepancy in surgical dosage is not merely a geometric artifact, as linear scaling models like Hirnschall et al. \cite{Hirnschall2022} fail to fully resolve, but a deterministic mechanical consequence of globe elongation that necessitates physics-based modeling.
\begin{table}[htbp]
\centering
\caption{Simulation Parameters for Sensitivity Analysis. The variations correspond to $\pm 1$ standard deviation (SD) derived from literature or estimated clinical variability}
\label{tab:sensitivity_params}
\begin{tabular}{l l c c c}
\toprule
\textbf{Parameter Class} & \textbf{Parameter Name} & \textbf{Nominal} & \textbf{\thead{Variation}} & \textbf{Ref.} \\
\midrule
\multirow{3}{*}{Geometric} & Eye Size (Radius) & 12.13~mm & $\pm 4\%$ (Iso-scaling) & \cite{Augusteyn2012} \\
 & LR Insertion (Ant.-Post.) & -6.5~mm & $\pm 0.8$ mm & \cite{apt1980anatomical} \\
 & LR Pulley (Vertical) & -0.3~mm & $\pm 2.0$ mm & \cite{kono2002} \\
\midrule
\multirow{3}{*}{Mechanical} & Tendon Slack & 8.4~mm & $\pm 5\%$ & \cite{Myers2014} \\
 & Passive Force & 1.0 (Normalized) & $\pm 15\%$ & \cite{Myers2014} \\
 & Max Force (Antagonist) & 1.0 (Normalized) & $\pm 15\%$ & \cite{Schutte2008} \\
\bottomrule
\end{tabular}
\end{table}

Our model clearly demonstrates a clinically significant difference of 0.9-1.5~mm in recession in high myopia compared with the emmetropic baseline, an essential difference for surgeons working with 0.5-mm increments. Indeed, as the axial length increases, muscle paths change, altering the moment arm of each muscle and explaining the difference in recession. \cite{Wright2009} showed that deviations as small as 2 PD can cause asthenopic symptoms and diplopia (incorrect image fusion). \cite{Hatt2012} made a comparative study of outcome criteria performance in adult strabismus surgery. They found that using only a motor criterion with an inferiority threshold of 10 PD in primary gaze is not enough. Overall, their study showed that applying motor criteria alone gives a 90 \% success rate, while combining motor criteria and diplopia gives a 67 \% success rate. Thus, the 8.23 PD under-correction may be acceptable for alignment but is not optimal for successful outcomes in highly myopic patients. Implemented in JAX, our model runs a complete optimized pipeline: from clinical measurement (PD), to patient-specific pathology ($\mathbf{F^*_{\text{max}}}$), to an optimized surgical plan ($\mathbf{p}^*_{\text{ins}}$).

However, modelling myopic eyes was limited due to a lack of data on the positions of the insertion and pulley points. Recent studies show the prevalence of pulleys in various pathologies involving strabismus \cite{Ludwig2025}, such as sagging or heavy eyes, and novel surgeries targeting those pulleys are emerging, necessitating an adaptation of our model for future plans. And while our sensitivity analysis demonstrates the theoretical robustness of the biomechanical mismatch, future work will focus on validating these predictions against post-operative data from myopic cohorts.

\section{Conclusion}\label{sec:conclusion}

We presented a differentiable, physics-based simulation of the eye for personalized planning of strabismus surgery. The model's differentiable nature allows for efficient gradient-based optimization to find the optimal surgical recession. 
The model was first validated against standard surgical guidelines for emmetropic eyes, demonstrating high consistency with the common surgical range. In cases where predictions slightly exceeded standard intervals, the mean deviation from the upper clinical boundary was only 0.15~mm, a clinically negligible difference compared to the standard 0.5~mm surgical increment.

Our results show that atypical anatomies, like high myopia, present biomechanical challenges not captured by standard tables. The model improves dosing accuracy for moderate deviations and demonstrates that the 10~mm safety cap is inadequate for large-angle strabismus in high myopes, often requiring bilateral surgery. This work represents a step toward safer, personalized planning that can both improve dose accuracy and prevent predictable surgical failures. Future research should enhance myopic eye modelling, commencing with MRI and OCT images from myopic patients.

\section*{Declarations}
\begin{itemize}
\item Funding: This project was funded by the State as part of France 2030 \includegraphics[height=3.0ex]{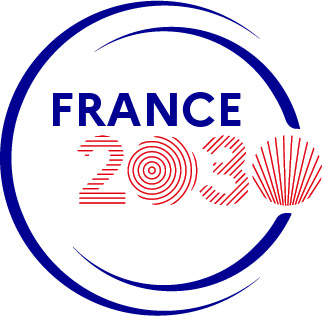}, and is part of the project PREMYOM. 
\item Competing interests: The authors have no conflicts of interest to declare
that are relevant to the content of this article.
\item Data availability: The data supporting this study are not publicly available.
\item Code availability: The code supporting this study is not publicly available.
\item Ethics approval: This article does not contain any studies with human participants or animals performed by any of the authors.
\item Informed Consent: Not applicable. This article does not contain any studies with
human participants performed by any of the authors.
\end{itemize}

\bigskip

\bibliographystyle{abbrv}
\bibliography{sn-bibliography}

\end{document}